\documentclass[a4paper,12pt]{article}
\usepackage[top=2.5cm,bottom=2.5cm,left=2.5cm,right=2.5cm]{geometry}
\usepackage{amsfonts}
\usepackage{latexsym,amsmath,amssymb,amsfonts,epsfig,psfrag,url,graphics,ifpdf,multicol}
\usepackage{color,colordvi,amscd,amsthm} 
\usepackage{tikz}
\usepackage[numbers,sort&compress]{natbib}
\usepackage{indentfirst,graphics,epsfig}
\usepackage{graphicx}
\usepackage{graphics}
\usepackage{graphicx,psfrag}
\usepackage{graphics}
\usepackage{mathrsfs}
\usepackage{tabularx}
\usepackage{booktabs}
\usepackage{floatrow}
\usepackage{multirow}
\usepackage{graphicx}
\usepackage{caption}
\usepackage[ruled]{algorithm2e}
\newfloatcommand{capbtabbox}{table}[][\FBwidth]

\newtheorem{theorem}{Theorem}
\newtheorem{theo}{Theorem}[section]
\newtheorem{lem}[theo]{Lemma}
\newtheorem{cl}{Claim}

 \theoremstyle{definition}

\theoremstyle{remark}

\def\pf{\noindent {\bf Proof.} }

\newcounter{casenum}[theo]

\newcounter{subcasenum}[theo]

\newcounter{claimnum}[theo]

\allowdisplaybreaks

\usepackage{indentfirst,subfig}
\begin{document}
    \thispagestyle{empty}
    \captionsetup[figure]{labelfont={bf},name={Fig.},labelsep=period}

\begin{center} {\large
On the maximum number of maximum dissociation sets in trees with given dissociation number}
\end{center}
\pagestyle{empty}

\begin{center}
{
  {\small Jianhua Tu$^1$, Lei Zhang$^2$, Junfeng Du$^{3,}$\footnote{Corresponding author.\\\indent \ \  E-mail: tujh81@163.com (J. Tu); 2018200896@mail.buct.edu.cn (L. Zhang); dujf1990@163.com (J. Du);}}\\[2mm]

{\small $^1$ School of Mathematics and Statistics, Beijing Technology and Business University, \\
	\hspace*{1pt} Beijing, P.R. China 100048} \\
{\small $^2$ School of Mathematics and Statistics, Beijing Institute of Technology, \\
	\hspace*{1pt} Beijing, P.R. China 100081} \\	
{\small $^3$ Department of mathematics, Beijing University of Chemical Technology, \\
\hspace*{1pt} Beijing, P.R. China 100029} \\[2mm]}

\end{center}

\begin{center}

\begin{abstract}
In a graph $G$, a subset of vertices is a dissociation set if it induces a subgraph with vertex degree at most 1. A maximum dissociation set is a dissociation set of maximum cardinality.
The dissociation number of $G$, denoted by $\psi(G)$, is the
cardinality of a maximum dissociation set of $G$. Extremal problems involving counting the number of a given type of substructure in a graph have been a hot topic of study in extremal graph theory throughout the last few decades. In this paper, we determine the maximum number of maximum dissociation sets in a tree with prescribed dissociation number and the extremal trees achieving this maximum value.


\vspace{5mm}

\noindent\textbf{Keywords:} Maximum dissociation set; Dissociation number; Tree; Extremal enumeration

\end{abstract}
\end{center}

\baselineskip=0.24in
\section{Introduction}

We consider only finite, simple, and undirected labeled graphs, and use standard terminology and notation as in \cite{Bondy2008}.

In a graph $G$, an \emph{independent set} is a set of vertices with no edges connecting them. 
A maximal independent set is an independent set that is not a proper subset of any other independent set. A maximum independent set is an independent set of maximum cardinality.
The \emph{independence number} of a graph $G$ is the cardinality of a maximum independent set of $G$ and is denoted by $\alpha(G)$.

In a graph $G$, a set of vertices is called a \emph{dissociation set} if it induces a subgraph with maximum degree at most 1. A \emph{maximum
dissociation set} of $G$ is a dissociation set of maximum cardinality. The \emph{dissociation number} of a graph $G$, denoted
by $\psi(G)$, is the cardinality of a maximum dissociation set of $G$. Clearly, an independent set is also a dissociation set, it follows that the concept of dissociation set is a natural generalization of the independent sets. The problem of finding a maximum dissociation set in a given graph was firstly studied by Yannakakis \cite{Yannakakis1981} and is NP-hard even in the class of bipartite and planar graphs. The complexity of the problem on some classes of graphs has been studied in
\cite{Alekseev2007,Bresar2011,Cameron2006,Orlovich2011}. A $k$-path vertex cover in a graph $G$ is a set of vertices intersecting every $k$-path in $G$, where a $k$-path is a not necessarily induced path of order $k$. It is easy to see that in a graph $G$, a set $F$ of vertices is a dissociation set if and only if its complement $V(G)\setminus F$ is a 3-path vertex cover. The problem of finding a 3-path vertex cover of minimum cardinality in a graph $G$ has received great attention throughout the past decade \cite{Bresar2011,Kardos2011,Katrenic2016,Tu2011,Xiao2017}.

In 1960s, Erd\H{o}s and Moser proposed the problem of determining the maximum number of \emph{maximal independent sets} in a general graph with prescribed number of vertices and extremal graphs achieving this maximum value. This problem was solved by Erd\H{o}s, Moon and Moser \cite{Moon1965}. Since then, there has been a great interest in a range of extremal problems involving counting the number of a given type of substructure in a graph. In particular, the extremal problems of determining the maximum number of maximal independent sets and maximum independent sets were extensively studied for various classes of graphs, including trees, forests, connected graphs, bipartite graphs, (connected) triangle-free graphs, (connected) graphs with at most $r$ cycles, unicyclic graphs etc. For these results, we refer to \cite{Chang1999, Furedi1987, Griggs1988, Jou2000, Koh2008, Liu1993, Sagan1988, Sagan2006, Wilf1986, Zito1991}. Researchers also studied the analogous extremal problems for some other graph substructures, including independent sets \cite{Sah2019,Zhao2010}, independent sets of a fixed size \cite{Gan2015, Chase2020}, minimal/minimum dominating sets \cite{Alvarado2019,Connolly2016}, minimal connected vertex covers \cite{Golovach2018}, minimal feedback vertex sets \cite{Couturier2012,Fomin2008}, maximal induced matchings \cite{Basavaraju2016}, etc. 

It should be pointed out that most of early work focused on the extremal problems on graphs with prescribed number of vertices. In 2019, Alvarado et al. \cite{Alvarado2019} studied the extremal problems on graphs where the independence number is prescribed instead of the number of vertices and proved that a tree $T$ of independence number $\alpha$ has at most $2^{\alpha-1}+1$ maximum independent sets. Mohr et al. \cite{Mohr2018, Mohr2020} studied the maximum number of maximum independent sets in general graphs, trees, subcubic trees and connected graphs with prescribed number of vertices and the independence number. Recently, Kirsch and Radcliffe \cite{Kirsch2021} determined the maximum number of cliques in a graph with prescribed number of edges and maximum degree.

Although dissociation sets were introduced and studied as early as 1980s, there is few results concerning the maximum number
of maximum dissociation sets and maximal dissociation sets in graphs. In \cite{Tu2021}, the first author, Zhang and Shi determined the maximum number of maximum dissociation sets in a tree with prescribed number of vertices and extremal trees achieving the maximum values. In this paper, we study the maximum number of maximum dissociation sets in a tree with prescribed dissociation number and prove that the maximum number of maximum dissociation sets in a tree of dissociation number $\psi$ is

\begin{center}
$\left\{
  \begin{array}{ll}
    1, &  \hbox{if $\psi=1$;} \\
    3^{\tfrac{\psi-1}{2}-1}+1, &  \hbox{if $\psi$ is odd and $\psi>1$;}\\
    3^{\tfrac{\psi}{2}-1}+\tfrac{\psi}{2}+1, & \hbox{if $\psi$ is even.}
  \end{array}
\right.$
\end{center}
We also give a complete characterization of extremal trees achieving the maximum values.


\section{The preparatory lemmas}

Let $G$ be a graph and $v$ a vertex in $G$. The neighborhood of $v$ is denoted by $N_G(v)$ and its cardinality denoted by $d_G(v)$. The closed neighborhood $N_G[v]$ of $v$ is $N_G(v)\cup\{v\}$. A vertex $u$ is called a leaf if $d_G(u)=1$. For a subset $X$ of vertices, the subgraph of $G$ induced on $X$ is denoted by $G[X]$. If $U$ is the set of vertices deleted, we write $G-U$ for the resulting subgraph.
We write $MD(G)$ and $\Phi(G)$ to denote the set of all maximum dissociation sets of $G$ and its cardinality, respectively. Let
\begin{align}
\Phi_v(G)&=|\{F\in MD(G): v\in F\}|, \notag\\
\Phi_{\overline{v}}(G)&=|\{F\in MD(G): v\notin F\}|, \notag\\
\Phi_v^0(G)&=|\{F\in MD(G): v\in F \text{ and } d_{G[F]}(v)=0\}|,\notag\\
\Phi_v^1(G)&=|\{F\in MD(G): v\in F \text{ and } d_{G[F]}(v)=1\}|,\notag
\end{align}
Clearly, $\Phi(G)=\Phi_v(G)+\Phi_{\overline{v}}(G)$ and
$\Phi_v(G)=\Phi_v^1(G)+\Phi_v^0(G)$.

Let $G$ and $H$ be two disjoint graphs, the union of $G$ and $H$, denoted by $G\cup H$, is the graph with vertex set $V(G\cup H)=V(G)\cup V(H)$ and edge set $E(G\cup H)=E(G)\cup E(H)$. We write $nG$ for the union of $n$ copies of disjoint graphs isomorphic to $G$.

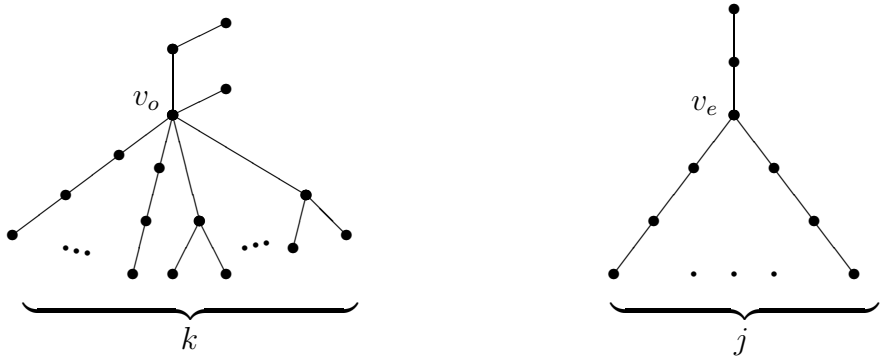
\begin{figure}[h]
    \begin{center}
        \begin{picture}(600,130)
        \multiput(120,80)(-5,-20){4}{\circle*{3.5}}
        \multiput(120,80)(-20,-15){4}{\circle*{3.5}}
        \multiput(120,80)(10,-40){2}{\circle*{3.5}}
        \multiput(120,80)(50,-30){2}{\circle*{3.5}}
        \multiput(130,40)(-10,-20){2}{\circle*{3.5}}
        \multiput(130,40)(10,-20){2}{\circle*{3.5}}
        \multiput(170,50)(-5,-20){2}{\circle*{3.5}}
        \multiput(170,50)(15,-15){2}{\circle*{3.5}}
        \multiput(120,80)(0,25){2}{\circle*{3.5}}
        \multiput(120,80)(20,10){2}{\circle*{3.5}}
        \multiput(120,105)(20,10){2}{\circle*{3.5}}
        \put(120,80){\line(-1,-4){15}}
        \put(120,80){\line(-4,-3){60}}
        \put(120,80){\line(5,-3){50}}
        \put(130,40){\line(-1,4){10}}
        \put(120,80){\line(2,1){20}}
        \put(120,105){\line(2,1){20}}
        \put(130,40){\line(-1,-2){10}}
        \put(130,40){\line(1,-2){10}}
        \put(170,50){\line(-1,-4){5}}
        \put(170,50){\line(1,-1){15}}
        \put(60,10){\ $\underbrace {\kern 125pt}$}
        \put(123,-8){$k$}
         \put(105,84){$v_o$}

        \put(120,80){\line(0,1){25}}
        \multiput(147,30)(4,1){3}{\circle*{2}}
        \multiput(80,30)(4,-1){3}{\circle*{2}}

        \multiput(330,80)(15,-20){4}{\circle*{3.5}}
        \multiput(330,80)(-15,-20){4}{\circle*{3.5}}
        \multiput(330,80)(0,20){3}{\circle*{3.5}}
        \put(330,80){\line(3,-4){45}}
        \put(330,80){\line(-3,-4){45}}
        \put(330,80){\line(0,1){40}}
        \put(280,10){\ $\underbrace {\kern 100pt}$}
        \put(330,-8){$j$}
        \multiput(315,20)(15,0){3}{\circle*{2}}
        \put(314,82){$v_e$}

        \caption{\label{1} A tree $T \in \mathcal{T}_1$ (left) and a tree $T \in \mathcal{T}_2$ (right).}
        \end{picture}
    \end{center}
\end{figure}

We define two special families of trees: $\mathcal{T}_1$ and $\mathcal{T}_2$. Every tree in the family $\mathcal{T}_1$ is obtained from $2P_2\cup kP_3$
by fixing a copy of $P_2$ and a vertex $v_o$ in it, and then adding edges joining $v_o$ to a vertex in each of the other copies of $P_2$ or $P_3$. Every tree in the family $\mathcal{T}_2$ is obtained from $(j+1)P_3$ by fixing a copy of $P_3$ and a leaf $v_e$ in it, and then adding edges joining $v_e$ to a leaf in each of the other copies of $P_3$. See Figure \ref{1}.

By a straightforward calculation, if $T\in \mathcal{T}_1$, then $\psi(T)=2k+3$, $\Phi(T)=3^k+1$ and $\Phi_{\overline{v_o}}(T)=3^k$; if $T\in \mathcal{T}_2$, then $\psi(T)=2j+2$, $\Phi(T)=3^j+j+2$ and $\Phi_{\overline{v_e}}(T)=3^j$.

Let
\begin{center}
$f(x)=\left\{
  \begin{array}{ll}
    3^{\tfrac{x-1}{2}-1}+1,&  \hbox{if $x$ is odd and $x>1$.}\\
    3^{\tfrac{x}{2}-1}+\tfrac{x}{2}+1,& \hbox{if $x$ is even;}
  \end{array}
\right.$
\end{center}
and
\begin{equation}
g(x)=\begin{cases}
     3^{\tfrac{x-1}{2}-1},&\text{if $x$ is odd},\\
     3^{\tfrac{x}{2}-1},&\text{if $x$ is even}.
     \end{cases}\notag
\end{equation}

The following two lemmas is straightforward to obtain. 

\begin{lem}
Let $T$ be a tree in the family $\mathcal{T}_1 $. Then

(1) $\psi(T)$ is an odd number;

(2) $\Phi(T)=f(\psi(T))$;

(3) $\Phi_{\overline{v_o}}(T)=g(\psi(T))$;

(4) and for any vertex $v$ other than $v_o$ in T, $\Phi_{\overline{v}}(T)< \Phi_{\overline{v_o}}(T).$
\end{lem}

\begin{lem}
Let $T$ be a tree in the family $\mathcal{T}_2$. Then

(1) $\psi(T)$ is an even number;

(2) $\Phi(T)=f(\psi(T))$;

(3) $\Phi_{\overline{v_e}}(T)=g(\psi(T))$;

(4) and for any vertex $v$ other than $v_e$ in $T$, $\Phi_{\overline{v}}(T)< \Phi_{\overline{v_e}}(T)$.
\end{lem}

\begin{lem}\label{lem1}
	Let $x$ and $m$ be two positive integers. If $m\ge 2$, then $f(x)<f(x+m)$.
\end{lem}
\pf When $m$ is even, $x$ and $x+m$ have the same parity
and this inequality $f(x)<f(x+m)$ obviously holds. Thus, we consider the case when $m$ is odd. It follows that $m\geq3$.

If $x$ is even, then

\[f(x+m)=3^{\tfrac{x+m-1}{2}-1}+1\ge 3^{\tfrac{x+2}{2}-1}+1 > 3^{\tfrac{x}{2}-1}+\tfrac{x}{2}+1.\]


If $x$ is odd, then

\[f(x+m)=3^{\tfrac{x+m}{2}-1}+\tfrac {x+m}{2}+1  \ge 3^{\tfrac{x+3}{2}-1}+\tfrac {x+3}{2}+1 > 3^{\tfrac{x-1}{2}-1}+1.\]

The proof is complete.\qed

\section{The main theorem}

Our main theorem is as follows.

\begin{theorem}\label{the1}
	Let $T$ be a tree of dissociation number $\psi$. Then
	
	(1) if $\psi = 1$, then $\Phi(T)=1$;
	
	(2) if $\psi$ is odd and $\psi\ge3$, then $\Phi(T)\le f(\psi)$ with equality if and only if $T$ is
	isomorphic to a tree in the family $\mathcal{T}_1$;
	
	(3) if $\psi$ is even, then $\Phi(T)\le f(\psi)$ with equality if and only if $T$ is
	isomorphic to a tree in the family $\mathcal{T}_2$.
	
\end{theorem}

\pf We call a tree {\it special} if it is isomorphic to a tree in the family $\mathcal{T}_1$ or $\mathcal{T}_2$. Suppose, to the contrary, that the theorem is false and there
exists some non special trees $T$ such that $\Phi(T)\geq f(\psi(T))$. Let $T_{0}$ be a tree which has
the smallest value of the dissociation number in all
of counterexamples. Let $\gamma:=\psi(T_0)$. Then $\Phi(T_{0})\geq f(\gamma)$ and $\gamma \geq 2$.\par

\begin{cl}\label{cl01}
Let $T$ be a tree of dissociation number $\psi$ and $v$ a vertex of $T$.
If $\psi<\gamma-1$ or $\psi=\gamma-1$ and $\psi$ is odd, then
$\Phi_{\overline{v}}(T)\le g(\psi)$.
\end{cl}

\noindent \textbf{Proof of Claim \ref{cl01}.} By the choice of $T_0$ and $\psi<\gamma$, we have that $\Phi(T)\leq f(\psi)$. We distinguish the following two
cases.

\textbf{Case 1.} $\psi$ is odd and $\psi\le\gamma-1$.\par

If $T$ is special, then $T\in \mathcal{T}_1$ and $\Phi_{\overline{v}}(T) \leq
\Phi_{\overline{v_o}}(T) = 3^{\tfrac{\psi-1}{2}-1}$;\par

If $T$ is not special, then $\Phi_{\overline{v}}(T)\le \Phi(T) <
3^{\tfrac{\psi-1}{2}-1} + 1$. Thus, $\Phi_{\overline{v}}(T)\leq
3^{\tfrac{\psi-1}{2}-1}$. \par

\textbf{Case 2.} $\psi$ is even and $\psi< \gamma-1$.\par

If $v$ is contained by all maximum dissociation sets of $T$, then
$\Phi_{\overline{v}}(T)=0$.

Now we assume that $\Phi_{\overline{v}}(T)>0$ and construct a new tree $T'$
from $T$ by attaching a new vertex $v'$ to $v$. Thus, $\psi(T')= \psi+1<\gamma$ and $\psi+1$ is odd. Furthermore, all maximum dissociation sets of $T'$ contain the vertex $v'$.

Since a maximum dissociation set in $T$  not containing $v$ can be extended to a maximum dissociation set
in $T'$ containing $v'$, $\Phi_{\overline{v}}(T) \leq \Phi_{v'}(T')=\Phi(T')$.

\begin{figure}[h]
    \begin{center}
        \begin{picture}(600,160)

        \multiput(120,100)(-20,-15){4}{\circle*{3.5}}
        \multiput(120,100)(-5,-20){4}{\circle*{3.5}}
        \multiput(120,100)(10,-40){2}{\circle*{3.5}}
        \multiput(120,100)(50,-30){2}{\circle*{3.5}}
        \multiput(130,60)(-10,-20){2}{\circle*{3.5}}
        \multiput(130,60)(10,-20){2}{\circle*{3.5}}
        \multiput(170,70)(-5,-20){2}{\circle*{3.5}}
        \multiput(170,70)(15,-15){2}{\circle*{3.5}}
        \multiput(120,100)(0,25){2}{\circle*{3.5}}
        \multiput(120,100)(20,10){2}{\circle*{3.5}}
        \multiput(120,125)(20,10){2}{\circle*{3.5}}
        \put(120,100){\line(-1,-4){15}}
        \put(120,100){\line(-4,-3){60}}
        \put(120,100){\line(5,-3){50}}
        \put(130,60){\line(-1,4){10}}
        \put(120,100){\line(2,1){20}}
        \put(120,125){\line(2,1){20}}
        \put(130,60){\line(-1,-2){10}}
        \put(130,60){\line(1,-2){10}}
        \put(170,70){\line(-1,-4){5}}
        \put(170,70){\line(1,-1){15}}
        \put(60,30){\ $\underbrace {\kern 125pt}$}
        \put(113,12){$\frac{\psi}{2}-1$}
        \put(138,142){$v'$}
        \put(109,128){$v$}
        \put(105,-5){$\Phi_{\overline {v}}(T)=1$}

        \put(120,100){\line(0,1){25}}
        \multiput(147,50)(4,1){3}{\circle*{2}}
        \multiput(80,50)(4,-1){3}{\circle*{2}}
        \multiput(340,100)(-5,-20){4}{\circle*{3.5}}
        \multiput(340,100)(-20,-15){4}{\circle*{3.5}}
        \multiput(340,100)(10,-40){2}{\circle*{3.5}}
        \multiput(340,100)(50,-30){2}{\circle*{3.5}}
        \multiput(350,60)(-10,-20){2}{\circle*{3.5}}
        \multiput(350,60)(10,-20){2}{\circle*{3.5}}
        \multiput(390,70)(-5,-20){2}{\circle*{3.5}}
        \multiput(390,70)(15,-15){2}{\circle*{3.5}}
        \multiput(340,100)(0,25){2}{\circle*{3.5}}
        \multiput(340,100)(20,10){2}{\circle*{3.5}}
        \multiput(340,125)(20,10){2}{\circle*{3.5}}
        \put(340,100){\line(-1,-4){15}}
        \put(340,100){\line(-4,-3){60}}
        \put(340,100){\line(5,-3){50}}
        \put(350,60){\line(-1,4){10}}
        \put(340,100){\line(2,1){20}}
        \put(340,125){\line(2,1){20}}
        \put(350,60){\line(-1,-2){10}}
        \put(350,60){\line(1,-2){10}}
        \put(390,70){\line(-1,-4){5}}
        \put(390,70){\line(1,-1){15}}
        \put(280,30){\ $\underbrace {\kern 125pt}$}
        \put(333,12){$\frac{\psi}{2}-1$}
        \put(360,115){$v'$}
        \put(329,104){$v$}
        \put(320,-5){$\Phi_{\overline {v}}(T)=3^{\frac{\psi}{2}-1}$}

        \put(340,100){\line(0,1){25}}
        \multiput(367,50)(4,1){3}{\circle*{2}}
        \multiput(300,50)(4,-1){3}{\circle*{2}}
        \caption{\label{2} Two possibilities for the vertex $v'$ by the structure of $T'$.}
\end{picture}
    \end{center}
\end{figure}
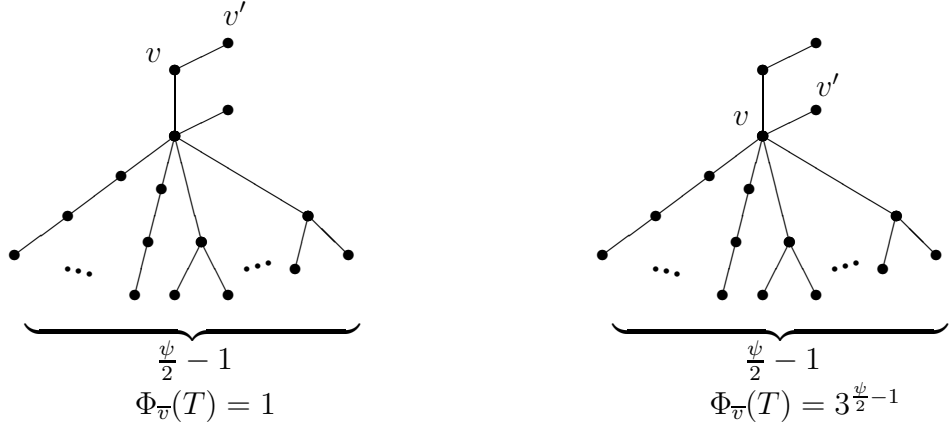

If $T'$ is not special, then $\Phi_{\overline{v}}(T)\le3^{\tfrac{\psi(T')-1}{2}-1}=3^{\tfrac{\psi}{2}-1}$.
If $T'$ is special and $\Phi_{v'}(T')= \Phi(T')=3^{\tfrac{\psi}{2}-1}+1$, then there are only two possibilities for $v'$ according to the structure of $T'$. See 
Figure \ref{2}. By direct calculation, in all two possibilities, we have that $\Phi_{\overline {v}}(T) \leq3^{\tfrac{\psi}{2}-1}$.

Consequently, we infer that if $\psi<\gamma-1$ or $\psi=\gamma-1$ and $\psi$ is odd, then
$\Phi_{\overline{v}}(T)\le g(\psi)$.\qed

\par

It can be checked by hand that the diameter of $T_0$ is at least 4, which implies that $\gamma\ge4$.

We change $T_0$ into a rooted tree by choosing an endvertex of a longest path in $T_0$ as the root. Recall that in a rooted tree $T$, the {\it level} of a vertex $v$ is the length of the path from the root to the vertex $v$. For a vertex $v$ other than the root, each vertex on the path from the root to the vertex $v$, not including the vertex $v$ itself, is called an ancestor of $v$; the parent of $v$ is the immediate ancestor of $v$. If $v$ is an ancestor of a vertex $u$, $u$ is called a descendant of $v$, the immediate descendants of $v$ are its children.


\begin{cl}\label{cl02}
Let $u$ be a non-leaf vertex of $T_{0}$. If all the children of $u$ are leaves, then $u$
has exactly one child.
\end{cl}

\noindent \textbf{Proof of Claim \ref{cl02}.} We assume that $u$ has $x$ children. Suppose,
to the contrary, that $x\ge 2$. Let $w$ be the parent of $u$. Let $T':=T_{0}-(N[u]\setminus \{w\})$. Now we consider the
following two cases.\par

\textbf{Case 1.} $x\ge 3$\par

In this case, $\psi(T')=\gamma-x$. Since a
maximum dissociation set in $T'$ can be extended in a unique
way to a maximum dissociation set in $T_{0}$, we have $\Phi(T_{0})=\Phi(T')$.  By Lemma \ref{lem1},
$\Phi(T_{0})=\Phi(T')\le f(\gamma-x) <f(\gamma)$, which contradicts the choice of $T_{0}$.\par

\textbf{Case 2.} $x= 2$\par

In this case, $\psi(T')=\gamma-2$. A maximum
dissociation set in $T'$ containing $w$ can be extended
in a unique way to a maximum dissociation set in $T_{0}$, and a
maximum dissociation set in $T'$ not containing $w$ can
be extended in three ways to a maximum dissociation set in $T_{0}$.
Thus,
\begin{equation}
\begin{aligned}
\Phi(T_{0})&= \Phi_{w}(T')+3 \cdot \Phi_{\overline{w}}(T')= \Phi(T')+2\cdot \Phi_{\overline{w}}(T'). \label{s1}
\end{aligned}
\end{equation}\par

If $\gamma$ is even, then by (\ref{s1}) and Claim \ref{cl01},
\begin{align}
\Phi(T_{0})&\le3^{\tfrac{\gamma-2}{2}-1}+\tfrac {\gamma-2}{2}+1+2\cdot 3^{\tfrac{\gamma-2}{2}-1}\notag=3^{\tfrac{\gamma}{2}-1}+\tfrac {\gamma}{2} <f(\gamma),\notag
\end{align}
which contradicts the choice of $T_{0}$.\par

If $\gamma$ is odd and $T'$ is not special, then by
(\ref{s1}) and Claim \ref{cl01}, $\Phi(T_{0})\le
3^{\tfrac{\gamma-3}{2}-1}+2\cdot
3^{\tfrac{\gamma-3}{2}-1}=3^{\tfrac{\gamma-1}{2}-1} <f(\gamma)$, which contradicts the choice of $T_{0}$.
Now we assume that $\gamma$ is odd and $T'$ is special. Since $T_{0}$ is not special, $w$ is other than vertex $v_o$ and
$\Phi_{\overline{w}}(T')<\Phi_{\overline{v_o}}(T')=3^{\tfrac{\gamma-3}{2}-1}$.
Thus, $\Phi(T_{0})< 3^{\tfrac{\gamma-3}{2}-1}+1+2\cdot
3^{\tfrac{\gamma-3}{2}-1}=3^{\tfrac{\gamma-1}{2}-1}+1 =f(\gamma)$,
which contradicts the choice of $T_{0}$.

In each case, we obtain a contradiction to the choice of $T_0$. The proof is complete.\qed

\vskip 3mm

Let $v$ be a leaf of maximum depth in $T_{0}$ and $vuwts$ a path in $T_0$.

\begin{cl} \label{cl03}
$u$ has exactly one child in $T_{0}$.
\end{cl}

\noindent \textbf{Proof of Claim \ref{cl03}.} The proof of this claim is directly obtained from Claim \ref{cl02} and the fact that $v$ is a leaf of maximum
depth in $T_{0}$. \qed

\begin{figure}[h]
    \begin{center}
        \begin{picture}(400,130)
         \multiput(200,70)(0,20){3}{\circle*{3.5}}
         \multiput(200,70)(-20,-30){2}{\circle*{3.5}}
         \multiput(200,70)(-50,-30){2}{\circle*{3.5}}
         \multiput(200,70)(20,-30){2}{\circle*{3.5}}
         \multiput(200,70)(50,-30){2}{\circle*{3.5}}
         \multiput(180,40)(0,-30){2}{\circle*{3.5}}
         \multiput(150,40)(0,-30){2}{\circle*{3.5}}
         \put(200,70){\line(0,1){55}}
         \put(200,70){\line(-2,-3){20}}
         \put(200,70){\line(-5,-3){50}}
         \put(200,70){\line(2,-3){20}}
         \put(200,70){\line(5,-3){50}}
         \put(180,40){\line(0,-1){30}}
         \put(150,40){\line(0,-1){30}}
         \multiput(160,25)(5,0){3}{\circle*{1.5}}
         \multiput(230,40)(5,0){3}{\circle*{1.5}}
         \put(140,47){\line(1,0){50}}
         \put(140,47){\line(0,-1){44}}
         \put(140,3){\line(1,0){50}}
         \put(190,47){\line(0,-1){44}}
         \put(210,47){\line(1,0){50}}
         \put(210,47){\line(0,-1){15}}
         \put(210,32){\line(1,0){50}}
         \put(260,47){\line(0,-1){15}}
         \put(140,38){$u$}
         \put(140,10){$v$}
         \put(160,-10){$P$}
         \put(230,20){$Q$}
         \put(183,68){$w$}
         \put(185,88){$t$}
         \put(185,108){$s$}

         \caption{\label{3} the vertex $w$ and its descendants}

         \end{picture}
    \end{center}
\end{figure}
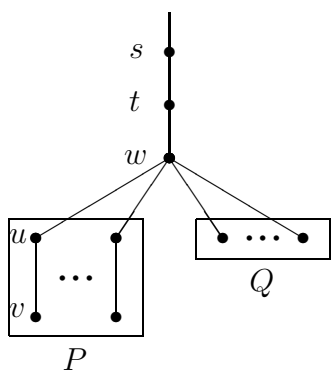

Assume that $w$ has $p$ children each of which is not a leaf and $q$ children
each of which is a leaf. Let $Q$ be the set of children of $w$ each of which is a leaf, and $P$ the set of all descendants of $w$ which are not
in $Q$. By Claim \ref{cl02}, each child of $w$ in $P$ has exactly one
child, in particular, $|P|=2p$. See Figure \ref{3}.

\begin{cl}\label{cl04}
$p=1$.
\end{cl}

\noindent \textbf{Proof of Claim \ref{cl04}.} Suppose, to the contrary, that
$p\ge2$. Let $T':=T_{0}-(P\cup Q\cup \{w\})$. Then
$\psi(T')=\gamma-2p-q\le \gamma-4$. Since a maximum dissociation
set in $T'$ can be extended in a unique way to a maximum
dissociation set in $T_{0}$, we have that $\Phi(T_{0})=\Phi(T')$. By Lemma \ref{lem1},
\begin{equation}
\Phi(T_{0})=\Phi(T')\le f(\psi(T')) <f(\gamma),\notag
\end{equation}
which contradicts the choice of $T_{0}$.\qed

\begin{cl}\label{cl05}
$q=0$.
\end{cl}

\noindent \textbf{Proof of Claim \ref{cl05}.} Suppose, to the contrary, that
$q\ge1$. Let $T':=T_{0}-(P\cup Q\cup \{w\})$. We consider
the following two cases.\par

\textbf{Case 1.} $q\ge 2$\par

In this case, $\psi(T')=\gamma-2-q\le \gamma-4$. Since a
maximum dissociation set in $T'$ can be extended in a unique
way to a maximum dissociation set in $T_{0}$, we have that $\Phi(T_{0})=\Phi(T')$. By Lemma \ref{lem1},
\begin{equation}
\Phi(T_{0})=\Phi(T')\le f(\psi(T')) <f(\gamma), \notag
\end{equation}
which contradicts the choice of $T_{0}$.\par

\textbf{Case 2.} $q=1$\par

In this case, $\psi(T')=\gamma-3$. Recall that $t$ is the parent of $w$. A maximum
dissociation set in $T'$ containing $t$ can be extended in
a unique way to a maximum dissociation set in $T_{0}$, and a maximum
dissociation set in $T'$ not containing $t$ can be extended in two ways to a maximum dissociation set in $T_{0}$. Thus,
\begin{equation}\label{s2}
\begin{aligned}
\Phi(T_0)=\Phi_t(T')+2 \cdot \Phi_{\overline{t}}(T')= \Phi(T')+ \Phi_{\overline{t}}(T').
\end{aligned}
\end{equation}


If $\gamma$ is even, then $\gamma-3$ is odd. By Claim \ref{cl01} and (\ref{s2}),
\begin{align}
\Phi(T_{0})\le3^{\tfrac{\gamma-4}{2}-1}+1+ 3^{\tfrac{\gamma-4}{2}-1}=2\cdot 3^{\tfrac{\gamma-4}{2}-1}+1\le 3^{\tfrac{\gamma}{2}-1} <f(\gamma),\notag
\end{align}
which contradicts the choice of $T_{0}$.\par

If $\gamma$ is odd, then $\gamma\geq5$ and $\gamma-3$ is even.
It is easy to check that for a tree $T$ of dissociation number 5,
$\Phi(T)\leq f(\gamma)=4$ with equality only if $T$ is special. Thus, we assume that $\gamma\ge7$. By Claim \ref{cl01} and
(\ref{s2}),
\begin{align}
\Phi(T_{0})\le3^{\tfrac{\gamma-3}{2}-1}+\tfrac{\gamma-3}{2}+1+ 3^{\tfrac{\gamma-3}{2}-1}\le 3^{\tfrac{\gamma-1}{2}-1} <f(\gamma),\notag
\end{align}
which contradicts the choice of $T_{0}$. \par

The proof is complete. \qed

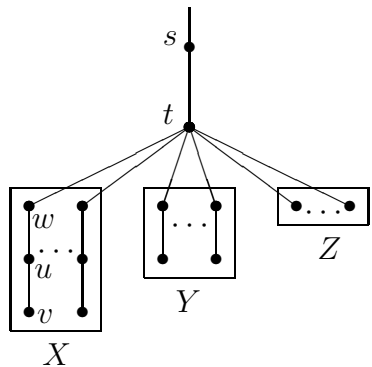
\begin{figure}[h]
    \begin{center}
        \begin{picture}(400,135)
         \multiput(200,80)(0,30){2}{\circle*{3.5}}
         \multiput(200,80)(-10,-30){2}{\circle*{3.5}}
         \multiput(200,80)(-40,-30){2}{\circle*{3.5}}
         \multiput(200,80)(-60,-30){2}{\circle*{3.5}}
         \multiput(200,80)(10,-30){2}{\circle*{3.5}}
         \multiput(200,80)(40,-30){2}{\circle*{3.5}}
         \multiput(200,80)(60,-30){2}{\circle*{3.5}}
         \multiput(190,50)(0,-20){2}{\circle*{3.5}}
         \multiput(210,50)(0,-20){2}{\circle*{3.5}}
         \multiput(160,50)(0,-20){3}{\circle*{3.5}}
         \multiput(140,50)(0,-20){3}{\circle*{3.5}}
         \put(200,80){\line(0,1){45}}
         \put(200,80){\line(2,-1){60}}
         \put(200,80){\line(4,-3){40}}
         \put(200,80){\line(1,-3){10}}
         \put(200,80){\line(-1,-3){10}}
         \put(200,80){\line(-4,-3){40}}
         \put(200,80){\line(-2,-1){60}}
         \put(190,50){\line(0,-1){20}}
         \put(210,50){\line(0,-1){20}}
         \put(160,50){\line(0,-1){40}}
         \put(140,50){\line(0,-1){40}}
         \put(133,57){\line(1,0){34}}
         \put(133,3){\line(1,0){34}}
         \put(133,57){\line(0,-1){54}}
         \put(167,57){\line(0,-1){54}}
         \put(183,57){\line(1,0){34}}
         \put(183,23){\line(1,0){34}}
         \put(183,57){\line(0,-1){34}}
         \put(217,57){\line(0,-1){34}}
         \put(233,57){\line(0,-1){14}}
         \put(267,57){\line(0,-1){14}}
         \put(233,57){\line(1,0){34}}
         \put(233,43){\line(1,0){34}}
         \put(145,-10){$X$}
         \put(195,10){$Y$}
         \put(193,40){$\cdots$}
         \put(248,30){$Z$}
          \put(243,45){$\cdots$}
         \put(190,81){$t$}
         \put(190,111){$s$}
         \put(141,42){$w$}
         \put(142,23){$u$}
         \put(143,30){$\cdots$}
         \put(143,6){$v$}
         \caption{\label{4} the vertex $t$ and its descendants}

         \end{picture}
    \end{center}
\end{figure}

Now, by Claim \ref{cl02}, Claim \ref{cl04} and Claim \ref{cl05}, we can assume that $t$ has $x$ children each of which has exactly two descendants, $y$ children each of which has exactly one descendant, and $z$ children each of which is a leaf. Let $X$, $Y$ and $Z$ be the set of the $x$, $y$ and $z$ children and all
descendants of them, respectively. In particular, $|X|=3x$ and $|Y|=2y$. (See Figure \ref{4})

\begin{cl}\label{cl06}
$z=0$,
\end{cl}

\noindent \textbf{Proof of Claim \ref{cl06}.} Suppose, to the contrary, that
$z\ge1$. Let $T'=T_{0}-X$. Then $\psi(T')=\gamma-2x$. Since $z\geq1$, if $F$ is a maximum dissociation set containing $t$ in $T'$, then $t$ must be an 1-degree vertex in $T'[F]$. Thus, a maximum dissociation
set in $T'$ containing $t$ can be extended in a unique way to a maximum dissociation set in $T_{0}$; a maximum dissociation
set in $T'$ not containing $t$ can be extended in $3^x$ ways
to a maximum dissociation set in $T_{0}$. So,
\begin{equation}\label{s3}
\begin{aligned}
\Phi(T_{0})=\Phi_t(T')+3^x \cdot \Phi_{\overline{t}}(T')= \Phi(T')+(3^x-1)\cdot \Phi_{\overline{t}}(T').
\end{aligned}
\end{equation}

Now, we consider the following two cases. \par

\textbf{Case 1.} $\gamma$ is even.\par

In this case, $\gamma-2x$ is even. By Claim \ref{cl01} and
(\ref{s3}),
\begin{align}
\Phi(T_{0})&\le3^{\tfrac{\gamma-2x}{2}-1}+\tfrac{\gamma-2x}{2}+1+ (3^x-1)\cdot 3^{\tfrac{\gamma-2x}{2}-1}\notag\\
&=3^{\tfrac{\gamma}{2}-1}+\tfrac{\gamma}{2}+1-x\le 3^{\tfrac{\gamma}{2}-1}+\tfrac{\gamma}{2} <f(\gamma),\notag
\end{align}
which contradicts the choice of $T_{0}$.\par

\textbf{Case 2.} $\gamma$ is odd.\par

In this case, $\gamma-2x$ is odd. If $T'$ is not special, then by
Claim \ref{cl01} and (\ref{s3}),
\begin{equation}
\Phi(T_{0})\le3^{\tfrac{\gamma-2x-1}{2}-1}+(3^x-1)\cdot3^{\tfrac{\gamma-2x-1}{2}-1}=3^{\tfrac{\gamma-1}{2}-1}
<f(\gamma), \notag
\end{equation}
which contradicts the choice of $T_{0}$.\par

Next, we assume that $T'$ is special. Since
$T_{0}$ is not special and $T'=T_0-X$, $t$ is other than the vertex
$v_o$ of $T'$ and $\Phi_{\overline{t}}(T')<\Phi_{\overline{v_o}}(T')=3^{\tfrac{\gamma-2x-1}{2}-1}$.
By Claim \ref{cl01} and (\ref{s3}),
\begin{equation}
\Phi(T_{0})<3^{\tfrac{\gamma-2x-1}{2}-1}+1+(3^x-1)\cdot
3^{\tfrac{\gamma-2x-1}{2}-1}=3^{\tfrac{\gamma-1}{2}-1}+1
=f(\gamma).\notag
\end{equation}
which contradicts the choice of $T_{0}$.\par

The proof is complete. \qed

\begin{cl}\label{cl07}
$y=0$.
\end{cl}

\noindent \textbf{Proof of Claim \ref{cl07}.} Suppose, to the contrary, that
$y\ge1$. We consider the following two cases. \par

\textbf{Case 1.} $y\ge 2$.\par

Let $T'=T_{0}-X$. Then $\psi(T')=\gamma-2x$. In this case, the vertex $t$ is in no maximum dissociation set of $T'$,
it follows that $\Phi(T')=\Phi_{\overline{t}}(T')$. Furthermore, a maximum dissociation set in
$T'$ can be extended in $3^x$ ways to a maximum dissociation set in
$T_{0}$.\par

When $\gamma$ is even, we have that $\Phi(T_{0})=3^x\cdot
\Phi_{\overline{t}}(T')\le 3^x\cdot
3^{\tfrac{\gamma-2x}{2}-1}=3^{\tfrac{\gamma}{2}-1}<f(\gamma)$, which contradicts the choice of $T_0$.

When $\gamma$ is odd, we have that $\Phi(T_{0})=3^x\cdot
\Phi_{\overline{t}}(T')\le 3^x\cdot
3^{\tfrac{\gamma-2x-1}{2}-1}=3^{\tfrac{\gamma-1}{2}-1}<f(\gamma)$, which contradicts the choice of $T_{0}$.\par

\textbf{Case 2.} $y=1$.\par

Let $T_1=T_{0}-(X\cup Y\cup \{t\})$, then
$\psi(T_1)=\gamma-2x-2$. Let $T_2=T_{0}-(X\cup Y)$, then $\psi(T_2)=\gamma-2x-2$ or $\psi(T_2)=\gamma-2x-1$. We consider the
following two subcases.\par

\textbf{Subcase 2.1.} $\psi(T_2)=\gamma-2x-2$\par

In this subcase, the vertex $t$ is in no maximum dissociation set of $T_0$. A maximum dissociation set in $T_2$ not containing $t$ can be extended in $3^x$ ways to a maximum dissociation set in $T_{0}$ and all maximum dissociation set in
$T_{0}$ are of such a form.

When $\gamma$ is even, by Claim \ref{cl01}, we have that $\Phi(T_{0})=3^x\cdot \Phi_{\overline{t}}(T_2)\le 3^x\cdot
3^{\tfrac{\gamma-2x-2}{2}-1}=3^{\tfrac{\gamma-2}{2}-1} <f(\gamma)$, which contradicts the choice of $T_0$.

When $\gamma$ is odd, we have that $\Phi(T_{0})=3^x\cdot \Phi_{\overline{t}}(T_2)\le 3^x\cdot
3^{\tfrac{\gamma-2x-3}{2}-1}=3^{\tfrac{\gamma-3}{2}-1} <f(\gamma)$, which contradicts the choice of $T_{0}$.

\textbf{Subcase 2.2.} $\psi(T_2)=\gamma-2x-1$\par

In this subcase, $\psi(T_2)=\psi(T_1)+1$. It follows that the vertex $t$ is in all maximum dissociation sets of $T_2$, and $\Phi(T_2)=\Phi_t(T_2)$. If $F$ is a maximum dissociation set in $T_2$ such that $t\in F$ and $d_{T_2[F]}(t)=0$, it can be
extended in $x+2$ ways to a maximum dissociation set in $T_{0}$; if $F$ is a
maximum dissociation set in $T_2$ such that $t\in F$ and $d_{T_2[F]}(t)=1$, it
can be extended in a unique way to a maximum dissociation set
in $T_{0}$, and all maximum dissociation sets of $T_{0}$ containing $t$ are of those forms. On the other hand, a maximum dissociation
set in $T_1$ can be extended in $3^x$ ways to a maximum dissociation
set in $T_{0}$ that does not contain $t$, and all maximum
dissociation sets of $T_{0}$ not containing $t$ are of that form. Thus,\\
\begin{equation}\label{s4}
\begin{aligned}
\Phi(T_{0})&=\Phi_{t}(T_0)+\Phi_{\overline{t}}(T_0)\\\
&=\{(x+2)\cdot\Phi_t^0(T_2)+\Phi_t^1(T_2)\}+3^x\cdot\Phi(T_1) \\\
 &= 3^x\cdot\Phi(T_1)+ \Phi(T_2)+ (x+1)\cdot\Phi_t^0(T_2)\\
 &= 3^x\cdot\Phi(T_1)+ \Phi(T_2)+ (x+1)\cdot\Phi_{\overline{s}}(T_2).
\end{aligned}
\end{equation}\par

When $\gamma$ is even, $\psi(T_1)=\gamma-2x-2$ is also even. If $\psi(T_1)=2$, then $|V(T_1)|\in\{2,3\}$ and
it is easy to see that $\Phi(T_0)\leq f(\gamma)=f(2x+4)=3^{x+1}+x+3$ with equality only if $T_0$ is special, which contradicts to the choice of $T_0$.
Thus $\psi(T_1)\ge 4$, $\gamma\ge 2x+6$ and $\gamma\geq8$. By Claim \ref{cl01} and (\ref{s4}),
\begin{align}
\Phi(T_{0})&\le 3^x\cdot (3^{\tfrac{\gamma-2x-2}{2}-1}+\tfrac{\gamma-2x-2}{2}+1)+3^{\tfrac{\gamma-2x-2}{2}-1}+1+(x+1)\cdot 3^{\tfrac{\gamma-2x-2}{2}-1}\notag\\
            &=3^{\tfrac{\gamma-4}{2}}+3^x\cdot(\tfrac{\gamma-2x}{2})+(x+2)\cdot3^{\tfrac{\gamma-2x-2}{2}-1}+1 \notag
\end{align}
Let $h_1(x)=3^x\cdot(\tfrac{\gamma-2x}{2})+(x+2)\cdot3^{\tfrac{\gamma-2x-2}{2}-1}+1$. Since $h_1''(x)>0$, $1\le x\le \frac{\gamma-6}{2}$ and $h_1(1)=h_1(\frac{\gamma-6}{2})$, we have $h_1(x)\le h_1(1)$ and
\begin{align}
\Phi(T_{0})\le 2\cdot 3^{\tfrac{\gamma-4}{2}}+ \tfrac{3\gamma}{2}-2\le3^{\tfrac{\gamma}{2}-1}+ \tfrac{\gamma}{2} <f(\gamma),\notag
\end{align}
which contradicts the choice of $T_{0}$.\par

When $\gamma$ is odd, $\psi(T_1)=\gamma-2x-2$ is also odd. If $\psi(T_1)=1$, then $|V(T_1)|=1$ and
$T_0$ is special, which contradicts to the choice of $T_0$. Thus, $\psi(T')\ge3$, $\gamma\ge 2x+5$ and $\gamma\geq7$. If $F$ is a
maximum dissociation set of $T_2$,  then $F-t$ is a maximum
dissociation set in $T_1$, which implies that $\Phi(T_2)\le\Phi(T_1)$.\par

Next we show that $\Phi(T_2)\le3^{\tfrac{\gamma-2x-3}{2}-1}$. If $T_1$ is not special, then $\Phi(T_2)\le\Phi(T_1)\le 3^{\tfrac{\gamma-2x-3}{2}-1}$. If $T_1$ is special, then $T_1$ is isomorphic to a tree in the family $\mathcal{T}_1$ and the fact that for every vertex $z$ of $T_1$ there exists a maximum dissociation set $F$ in $T_1$ such that $z\in F$ and $d_{T_1[F]}(z)=1$ implies that $\Phi(T_2)< \Phi(T_1)$. In both cases, we have that $\Phi(T_2)\le3^{\tfrac{\gamma-2x-3}{2}-1}$.
By Claim \ref{cl01} and (\ref{s4}),
\begin{align}
\Phi(T_{0})&\le 3^x\cdot (3^{\tfrac{\gamma-2x-3}{2}-1}+1)+3^{\tfrac{\gamma-2x-3}{2}-1}+(x+1)\cdot 3^{\tfrac{\gamma-2x-3}{2}-1}\notag\\
            &=3^{\tfrac{\gamma-5}{2}}+3^x+(x+2)\cdot3^{\tfrac{\gamma-2x-3}{2}-1} \notag
\end{align}
Let $h_2(x)=3^x+(x+2)\cdot3^{\tfrac{\gamma-2x-3}{2}-1}$. Since $h_2''(x)>0$, $1\le x\le \frac{\gamma-5}{2}$ and $h_2(\frac{\gamma-5}{2})\ge h_2(1)$, we have that $h_2(x)\leq h_2(\frac{\gamma-5}{2})$ and
\begin{align}
\Phi(T_{0})\le 2\cdot 3^{\tfrac{\gamma-5}{2}}+\tfrac{\gamma-1}{2}\le3^{\tfrac{\gamma-1}{2}-1} <f(\gamma),\notag
\end{align}
which contradicts the choice of $T_{0}$.\par

In each case, we obtain a contradiction to the choice of $T_{0}$. We complete the proof of Claim \ref{cl07}. \qed

We are now in a position to derive a final contradiction and
consider the following two cases. \par

\textbf{Case 1.} $\gamma$ is even.\par

Let $T_a=T_{0}-X$ and $T_b=T_{0}-(X\cup \{t\})$. Then
$\psi(T_a)=\gamma-2x$, and $\psi(T_b)=\gamma-2x-1$ or $\psi(T_b)=\gamma-2x$. If $\psi(T_a)=2$, then $|V(T_a)|\in\{2,3\}$ and
it is easy to see that $\Phi(T_0)\leq f(\gamma)=f(2x+2)=3^{x}+x+2$ with equality only if $T_0$ is special, which contradicts to the choice of $T_0$.
Thus $\psi(T_a)\ge 4$, $\gamma\ge 2x+4$ and $\gamma\geq6$. Now, we consider the following two subcases.\par

\textbf{Subcase 2.1.} $\psi(T_b)=\gamma-2x-1$\par

In this subcase, the vertex $t$ is in all maximum dissociation sets of $T_a$ and $\Phi(T_a)=\Phi_t(T_a)$.
If $F$ is a maximum dissociation set in $T_a$ such that $d_{T_a[F]}(t)=0$, it can be extended in $x+1$
ways to a maximum dissociation set in $T_{0}$; if $F$ is a maximum
dissociation set in $T_a$ such that $d_{T_a[F]}(t)=1$, it can be
extended in a unique way to a maximum dissociation set in $T_{0}$.
Since $\Phi_t^0(T_a) = \Phi_{\overline{s}}(T_a)\le
3^{\tfrac{\gamma-2x}{2}-1}$, we have
\begin{align}
\Phi(T_{0})&=\Phi_t^1(T_a)+(x+1)\cdot\Phi_t^0(T_a)\notag\\
&=\Phi(T_a)+x\cdot \Phi_t^0(T_a)\le3^{\tfrac{\gamma-2x}{2}-1}+ \tfrac{\gamma-2x}{2}+1+x\cdot 3^{\tfrac{\gamma-2x}{2}-1}\notag
\end{align}
Let $\ell_1(x)=3^{\tfrac{\gamma-2x}{2}-1}+ \tfrac{\gamma-2x}{2}+1+x\cdot 3^{\tfrac{\gamma-2x}{2}-1}$. Since $\ell_1''(x)>0$, $1\leq x\leq \frac{\gamma-4}{2}$ and $\ell_1(1)\geq \ell_1(\frac{\gamma-4}{2})$, we have
\begin{align}
\Phi(T_{0})\le 2\cdot 3^{\tfrac{\gamma-4}{2}}+\tfrac{\gamma}{2}\le 3^{\tfrac{\gamma}{2}-1} <f(\gamma), \notag
\end{align}
which contradicts the choice of $T_{0}$.\par

\textbf{Subcase 2.2.} $\psi(T_b)=\gamma-2x$\par

In this subcase, $\psi(T_a)=\psi(T_b)$. Thus, for every maximum dissociation set $F$ in $T_b$, $s\in F$ and $d_{T_b[F]}(s)=1$.
Let $N_{T_b}(s)=\{t_1,\cdots,t_k\}$. Next, we will prove that there exists a vertex in $N_{T_b}(s)$ such that it is in all maximum dissociation sets of $T_b$.

Without loss of generality, suppose to the contrary that there are two maximum dissociation sets $F_1$ and $F_2$ in $T_b$ containing $t_1$ and $t_2$, respectively. Let $D_1,\cdots,D_k$ be the connected components of $T_b-s$ such that $t_i\in V(D_i)$. It is easy to see that $F_1\cap V(D_1)$ is a maximum dissociation set of $D_1$ and every maximum dissociation set of $D_1$ contains the vertex $t_1$. Since $F_2\cap V(D_1)$ is a dissociation set of $D_1$ not containing $t_1$, $|F_2\cap V(D_1)|< |F_1\cap V(D_1)|$. Now, $[(F_2\setminus (F_2\cap V(D_1)))\cup (F_1\cap V(D_1))]\setminus\{s\}$ is a maximum dissociation set in $T_b$ that does not contain $s$, this is a contradiction. Thus, we verify that there exists a vertex in $N_{T_b}(s)$, say $t_1$, such that it is in all maximum dissociation sets of $T_b$.

If $F$ is a maximum dissociation set in $T_b$, then $F_1=(F\cup\{t\})\setminus\{s\}$ is a maximum dissociation set of $T_a$ such that $d_{T_a[F_1]}(t)=0$, and $F_2=(F\cup\{t\})\setminus\{t_1\}$ is a maximum dissociation set in $T_a$ such that $d_{T_a[F_2]}(t)=1$. On the other hand, $F$ is also a maximum dissociation set in $T_a$ that does not contain $t$. Thus, we have
\begin{align}
\Phi_{\overline{t}}(T_a)=\Phi(T_b)\leq \min\{\Phi_t^0(T_a),\Phi_t^1(T_a)\}\notag,
\end{align}
Since $\Phi(T_a)=\Phi_{\overline{t}}(T_a)+\Phi_t^0(T_a)+\Phi_t^1(T_a)\le
3^{\tfrac{\gamma-2x}{2}-1}+\tfrac{\gamma-2x}{2}+1$, we have
\begin{align}
\Phi_{\overline{t}}(T_a)&\le \frac{1}{3}(
3^{\tfrac{\gamma-2x}{2}-1}+\tfrac{\gamma-2x}{2}+1),\label{f}\\
2\Phi_{\overline{t}}(T_a)+\Phi_t^0(T_a)&\le\Phi(T_a)\leq 3^{\tfrac{\gamma-2x}{2}-1}+\tfrac{\gamma-2x}{2}+1\label{g}.
\end{align}
A maximum dissociation set in $T_a$ not containing $t$ can be extended in $3^x$ ways to a maximum
dissociation set in $T_{0}$, a maximum dissociation set $F$ in $T_a$
that contains $t$ such that $d_{T_a[F]}(t)=0$ can be extended in
$x+1$ ways to a maximum dissociation set in $T_{0}$, and a maximum
dissociation set $F$ in $T_a$ that contains $t$ such that
$d_{T'[F]}(t)=1$ can be extended in a unique way to a maximum
dissociation set in $T_{0}$. Since all maximum dissociation sets of
$T_{0}$ are of such forms, we obtain
\begin{align}
\Phi(T_{0})&=3^x\cdot\Phi_{\overline{t}}(T_a)+(x+1)\cdot \Phi_t^0(T_a)+\Phi_t^1(T_a) \notag \\
              &=\Phi(T_a)+(3^x-1)\cdot \Phi_{\overline{t}}(T_a)+x\cdot \Phi_t^0(T_a). \notag
\end{align}
By (\ref{f}) and (\ref{g}), consider the following linear programming:
\begin{align}
\max.\ &(3^x-1)\cdot \Phi_{\overline{t}}(T_a)+x\cdot \Phi_t^0(T_a) \notag \\
\text{s. t.} \ \ &\Phi_{\overline{t}}(T_a) \le \tfrac{1}{3}\cdot (3^{\tfrac{\gamma-2x}{2}-1}+\tfrac{\gamma-2x}{2}+1) \notag \\
&\Phi_t^0(T_a)+2\Phi_{\overline{t}}(T_a)\le 3^{\tfrac{\gamma-2x}{2}-1}+\tfrac{\gamma-2x}{2}+1 \notag \\
&\Phi_t^0(T_a)\ge0,\ \Phi_{\overline{t}}(T_a) \ge 0\notag
\end{align}
The linear programming has an optimal solution $(\Phi_{\overline{t}}(T_a),\Phi_t^0(T_a))=(\tfrac{1}{3} (3^{\tfrac{\gamma-2x}{2}-1}+\tfrac{\gamma-2x}{2}+1), \tfrac{1}{3}(3^{\tfrac{\gamma-2x}{2}-1}+\tfrac{\gamma-2x}{2}+1))$. Thus,
\begin{align}
\Phi(T_{0})&=\Phi(T_a)+(3^x-1)\cdot \Phi_{\overline{t}}(T_a)+x\cdot \Phi_t^0(T_a) \notag\\
&\le 3^{\tfrac{\gamma-2x}{2}-1}+\tfrac{\gamma-2x}{2}+1+(3^x-1+x)\cdot  \tfrac{1}{3}\cdot (3^{\tfrac{\gamma-2x}{2}-1}+\tfrac{\gamma-2x}{2}+1)\notag\\
              &=(3^{x-1}+\tfrac{x+2}{3})(3^{\tfrac{\gamma-2x}{2}-1}+\tfrac{\gamma-2x}{2}+1)\notag\\
              &=3^{\tfrac{\gamma}{2}-2}+3^{x-1}\cdot (\tfrac{\gamma-2x+2}{2})+(\tfrac{x+2}{3})\cdot (3^{\tfrac{\gamma-2x}{2}-1}+\tfrac{\gamma-2x}{2}+1) \notag
\end{align}
Let $\ell_2(x)=3^{x-1}\cdot (\tfrac{\gamma-2x+2}{2})+(\tfrac{x+2}{3})\cdot (3^{\tfrac{\gamma-2x}{2}-1}+\tfrac{\gamma-2x}{2}+1)$. Since $\ell_2''(x)>0$, $1\leq x\leq \frac{\gamma-4}{2}$ and $\ell_2(1)=\ell_2(\frac{\gamma-4}{2})$, we have
\begin{align}
\Phi(T_{0})\le2\cdot 3^{\tfrac{\gamma}{2}-2}+\gamma\le 3^{\tfrac{\gamma}{2}-1}+\tfrac{\gamma}{2} <f(\gamma),\notag
\end{align}
which contradicts the choice of $T_{0}$.\par

\textbf{Case 2.} $\gamma$ is odd.\par

Let $T_a=T_{0}-X$. Clearly, $\psi(T_a)=\gamma-2x\ge 3$. Similarly, we have
\begin{align}
\Phi(T_{0})&=3^x\cdot \Phi_{\overline{t}}(T_a)+(x+1)\cdot \Phi_t^0(T_a)+\Phi_t^1(T_a)\notag
\end{align}

If $T_a$ is not special, then
\begin{align}
    \Phi(T_{0})&\le3^x\cdot[\Phi_{\overline{t}}(T_a)+\Phi_t^0(T_a)+\Phi_t^1(T_a)]\notag \\
 &=3^x \cdot \Phi(T_a)\le 3^x\cdot 3^{\tfrac{\gamma-2x-1}{2}-1}=3^{\tfrac{\gamma-1}{2}-1} <f(\gamma).\notag
\end{align}
%
which contradicts the choice of $T_{0}$.\par

\begin{figure}[h]
    \begin{center}
        \begin{picture}(360,140)
        \multiput(120,90)(-5,-20){4}{\circle*{3.5}}
        \multiput(120,90)(-20,-15){4}{\circle*{3.5}}
        \multiput(120,90)(10,-40){2}{\circle*{3.5}}
        \multiput(120,90)(50,-30){2}{\circle*{3.5}}
        \multiput(130,50)(-10,-20){2}{\circle*{3.5}}
        \multiput(130,50)(10,-20){2}{\circle*{3.5}}
        \multiput(170,60)(-5,-20){2}{\circle*{3.5}}
        \multiput(170,60)(15,-15){2}{\circle*{3.5}}
        \multiput(120,90)(0,25){2}{\circle*{3.5}}
        \multiput(120,90)(20,10){2}{\circle*{3.5}}
        \multiput(120,115)(20,10){2}{\circle*{3.5}}
        \put(120,90){\line(-1,-4){15}}
        \put(120,90){\line(-4,-3){60}}
        \put(120,90){\line(5,-3){50}}
        \put(130,50){\line(-1,4){10}}
        \put(120,90){\line(2,1){20}}
        \put(120,115){\line(2,1){20}}
        \put(130,50){\line(-1,-2){10}}
        \put(130,50){\line(1,-2){10}}
        \put(170,60){\line(-1,-4){5}}
        \put(170,60){\line(1,-1){15}}
        \put(60,20){\ $\underbrace {\kern 125pt}$}
        \put(100,2){$\frac{\gamma-2x-1}{2}-1$}
        \put(138,132){$t^1$}
        \put(138,106){$t^2$}
        \put(182,34){$t^3$}
        \put(53,32){$t^4$}

        \put(120,90){\line(0,1){25}}
        \multiput(147,40)(4,1){3}{\circle*{2}}
        \multiput(80,40)(4,-1){3}{\circle*{2}}
        \multiput(240,100)(0,-30){3}{\circle*{3.5}}
        \multiput(265,100)(0,-30){3}{\circle*{3.5}}
        \multiput(290,100)(0,-30){3}{\circle*{3.5}}
        \put(240,100){\line(0,-1){60}}
        \put(265,100){\line(0,-1){60}}
        \put(290,100){\line(0,-1){60}}
        \multiput(273,70)(5,0){3}{\circle*{2}}
        \put(233,107){\line(1,0){64}}
        \put(233,33){\line(1,0){64}}
        \put(233,33){\line(0,1){74}}
        \put(297,33){\line(0,1){74}}
        \put(265,18){$X$}
        \put(140,125){\line(4,-1){100}}
        \put(140,125){\line(6,-1){150}}
        \put(140,125){\line(5,-1){125}}
        \caption{\label{5} Four possibilities for the structure of $T$}
\end{picture}
    \end{center}
\end{figure}
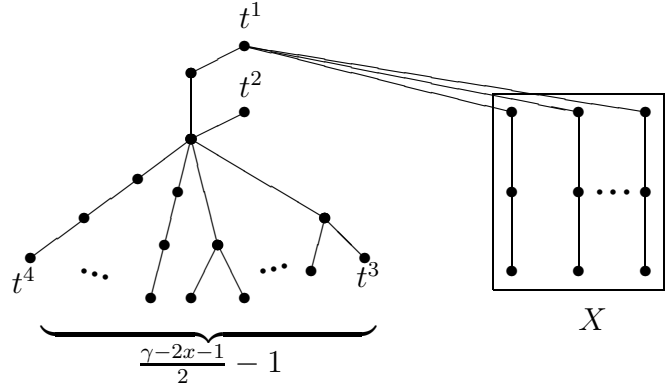

If $T_a$ is special, then there are four possibilities for the vertex $t$. See Figure \ref{5}.
We use $T^1$, $T^2$, $T^3$ and $T^4$ to
represent the four possibilities of structures of $T_0$ such that $t=t^i$ in $T^i$. It should be noted that when $T_0\cong T^3 \text{\ or\ } T^4$, $\psi(T_0)=\gamma\geq7$.

By direct calculation, we obtain
\begin{align}
&\Phi(T^1)=3^{\tfrac{\gamma-2x-1}{2}-1}+x+1 \le3^{\tfrac{\gamma-1}{2}-1} <f(\gamma),\notag\\
&\Phi(T^2)=(x+1)\cdot 3^{\tfrac{\gamma-2x-1}{2}-1}+1 \le3^{\tfrac{\gamma-1}{2}-1} <f(\gamma),\notag\\
&\Phi(T^3)= 3^{\tfrac{\gamma-5}{2}}+(x+2)\cdot 3^{\tfrac{\gamma-2x-1}{2}-2}+x+1 \le 3^{\tfrac{\gamma-1}{2}-1} <f(\gamma),\notag\\
&\Phi(T^4)=3^{\tfrac{\gamma-5}{2}}+(x+2)\cdot
3^{\tfrac{\gamma-2x-1}{2}-2}+1\le 3^{\tfrac{\gamma-1}{2}-1}
<f(\gamma).\notag
\end{align}
All of these contradict the choice of $T_{0}$.\par

In each case, we obtain a contradiction to the choice of $T_{0}$. The proof of Theorem \ref{the1} is complete.\qed

\section*{Acknowledgments}
\noindent The work is supported by Research Foundation for Advanced Talents of Beijing Technology and Business University.


\end{document}